\documentclass[12pt]{article}

\usepackage{fullpage}
\usepackage{amssymb}
\usepackage{amsmath}
\usepackage{epsfig}
\usepackage{algorithm}
\usepackage{algorithmic}
\usepackage{latexsym}
\usepackage{amsmath}
\usepackage{amsfonts}
\usepackage{graphicx}
\usepackage[normalem]{ulem}
\usepackage{array}
\usepackage{color}\usepackage[dvipsnames]{xcolor}
\usepackage{tikz}\usetikzlibrary{patterns}
\usepackage{array,multirow}
\usepackage[
  bookmarks=false, 
  colorlinks,
  citecolor=brown!70!black,
  linkcolor=brown!80!black,
  urlcolor=blue!70!black,
]{hyperref}
\usepackage{relsize,exscale}

\newtheorem{lem}{Lemma}

\newtheorem{thm}{Theorem}

\newtheorem{fac}{Fact}

\newcommand{\comm}[1]{}

\newcolumntype{C}[1]{>{\centering\arraybackslash }b{#1}}

\def\A{0.5cm}

 \def\C{2.5cm}
 
 \def\E{4.5cm}

 \def\G{6.5cm}
 
 \def\I{8.5cm}

 \def\K{10.5cm}
 
\def\M{12.5cm}
 
 \def\O{14.5cm}\def\Q{16.5cm}\def\S{18.5cm}
\def\U{20.5cm}\def\W{22.5cm}\def\Y{24.5cm}\def\ZZ{26.5cm}\def\Za{28.5cm}\def\Zb{30.5cm}\def\Zc{32.5cm}\def\Zd{34.5cm}\def\Ze{36.5cm}

\title{Enumeration of \L{}ukasiewicz paths modulo some patterns}
\author{\large Jean-Luc {Baril}, Sergey Kirgizov and Armen Petrossian\\  LE2I, Universit\'e de Bourgogne\\
         B.P. 47 870, 21078 DIJON-Cedex France      \\
        {\tt e-mail: \{barjl, armen.petrossian\}@u-bourgogne.fr},\\{\tt kerzolster@gmail.com}}

\begin{document}
\maketitle

\begin{abstract}  For any pattern $\alpha$ of length at most two,  we enumerate equivalence classes of \L{}ukasiewicz paths of length $n\geq 0$ where two paths are equivalent whenever the
occurrence positions of $\alpha$ are identical on these paths. As a byproduct, we give a constructive bijection between Motzkin paths and  some equivalence classes of  \L{}ukasiewicz paths.
\end{abstract}


\section{Introduction and notations}
In the literature, lattice paths are widely studied. Their enumeration is a very active field in combinatorics, and they have many applications in other research domains as computer science,
 biology and physics \cite{Moh,Nar}.  Dyck and Motzkin paths are the most often considered. This is partly due to the fact that they are respectively counted by the famous Catalan
 and Motzkin numbers (see \href{https://oeis.org/A000108}{A000108}  and \href{https://oeis.org/A001006}{A001006} in the on-line encyclopedia of integer sequences \cite{Sloa}). Almost always, these paths are enumerated according to several
  parameters and statistics (see for instance \cite{Deu, Man1,Man2,Mer,Pan,Pea,Sap,Sap1,Sun} for Dyck paths and \cite{Barc,Bre,Don,Dra,Man4,Pro,Sap2} for Motzkin paths). Also, many
  one-to-one correspondences  have been found between lattice paths and some  combinatorial objects such as Young tableaux, pattern avoiding permutations, bargraphs, RNA shapes and so on
 \cite{Sta}. Recently a new approach has been introduced for studying statistics on lattice paths. It consists in determining the cardinality of the quotient set generated by an equivalence
   relation based on the positions of a given pattern: {\it two paths belong to the same equivalence class whenever the positions of  occurrences of a given pattern are identical on these paths.}
   Enumerating results are provided for the quotient sets of Dyck, Motzkin and Ballot paths for  patterns of length at most three (see respectively \cite{Bar1}, \cite{Bar2} and \cite{Man}). The purpose
   of this present paper is to extend these studies for  \L{}ukasiewicz paths that  naturally  generalizes Dyck and Motzkin paths. As a byproduct, we show how Motzkin paths are in one-to-one
   correspondence with some equivalence classes of   \L{}ukasiewicz paths.

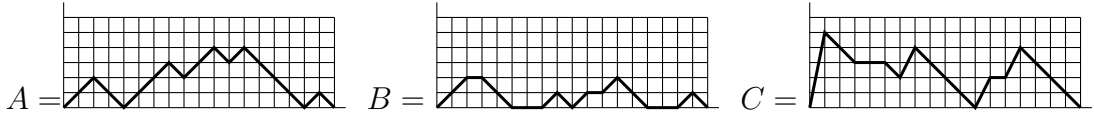
\begin{figure}[h]
 \begin{center}
 $A=$\begin{tikzpicture}[scale=0.1]
            \draw (\A,\A)-- (38,\A);
             \draw (\A,\E)-- (\Ze,\E);
              \draw (\A,\C)-- (\Ze,\C);
               \draw (\A,\G)-- (\Ze,\G);
               \draw (\A,\I)-- (\Ze,\I);
              \draw (\A,\K)-- (\Ze,\K);
               \draw (\A,\M)-- (\Ze,\M);
            \draw (\A,\A) -- (\A,\O);
             \draw (\C,\A) -- (\C,\M);\draw (\E,\A) -- (\E,\M);\draw (\G,\A) -- (\G,\M);
             \draw (\I,\A) -- (\I,\M);\draw (\K,\A) -- (\K,\M);\draw (\M,\A) -- (\M,\M);
             \draw (\O,\A) -- (\O,\M);\draw (\Q,\A) -- (\Q,\M);\draw (\S,\A) -- (\S,\M);
             \draw (\U,\A) -- (\U,\M);\draw (\W,\A) -- (\W,\M);\draw (\Y,\A) -- (\Y,\M);
             \draw (\ZZ,\A) -- (\ZZ,\M);
             \draw (\Za,\A) -- (\Za,\M);
             \draw (\Zb,\A) -- (\Zb,\M);
             \draw (\Zc,\A) -- (\Zc,\M);
             \draw (\Zd,\A) -- (\Zd,\M);
             \draw (\Ze,\A) -- (\Ze,\M);
            \draw[solid,line width=0.4mm] (\A,\A)--(\C,\C) -- (\E,\E) -- (\G,\C) -- (\I,\A) -- (\O,\G) -- (\Q,\E) -- (\U,\I) -- (\W,\G) --(\Y,\I) --  (\Zc,\A) -- (\Zd,\C) -- (\Ze, \A);
         \end{tikzpicture}
     $~B=$    \begin{tikzpicture}[scale=0.1]
            \draw (\A,\A)-- (38,\A);
             \draw (\A,\E)-- (\Ze,\E);
              \draw (\A,\C)-- (\Ze,\C);
               \draw (\A,\G)-- (\Ze,\G);
               \draw (\A,\I)-- (\Ze,\I);
              \draw (\A,\K)-- (\Ze,\K);
               \draw (\A,\M)-- (\Ze,\M);
            \draw (\A,\A) -- (\A,\O);
             \draw (\C,\A) -- (\C,\M);\draw (\E,\A) -- (\E,\M);\draw (\G,\A) -- (\G,\M);
             \draw (\I,\A) -- (\I,\M);\draw (\K,\A) -- (\K,\M);\draw (\M,\A) -- (\M,\M);
             \draw (\O,\A) -- (\O,\M);\draw (\Q,\A) -- (\Q,\M);\draw (\S,\A) -- (\S,\M);
             \draw (\U,\A) -- (\U,\M);\draw (\W,\A) -- (\W,\M);\draw (\Y,\A) -- (\Y,\M);
             \draw (\ZZ,\A) -- (\ZZ,\M);
             \draw (\Za,\A) -- (\Za,\M);
             \draw (\Zb,\A) -- (\Zb,\M);
             \draw (\Zc,\A) -- (\Zc,\M);
             \draw (\Zd,\A) -- (\Zd,\M);
             \draw (\Ze,\A) -- (\Ze,\M);
            \draw[solid,line width=0.4mm] (\A,\A)--(\C,\C) -- (\E,\E) -- (\G,\E) -- (\I,\C) -- (\K,\A) -- (\M,\A) -- (\O,\A) -- (\Q,\C) -- (\S,\A) -- (\U,\C)  -- (\W,\C) --(\Y,\E) -- (\ZZ,\C) --(\Za,\A) -- (\Zb,\A) -- (\Zc,\A) -- (\Zd,\C) -- (\Ze, \A);
         \end{tikzpicture}
          $~C=$    \begin{tikzpicture}[scale=0.1]
            \draw (\A,\A)-- (38,\A);
             \draw (\A,\E)-- (\Ze,\E);
              \draw (\A,\C)-- (\Ze,\C);
               \draw (\A,\G)-- (\Ze,\G);
               \draw (\A,\I)-- (\Ze,\I);
              \draw (\A,\K)-- (\Ze,\K);
               \draw (\A,\M)-- (\Ze,\M);
            \draw (\A,\A) -- (\A,\O);
             \draw (\C,\A) -- (\C,\M);\draw (\E,\A) -- (\E,\M);\draw (\G,\A) -- (\G,\M);
             \draw (\I,\A) -- (\I,\M);\draw (\K,\A) -- (\K,\M);\draw (\M,\A) -- (\M,\M);
             \draw (\O,\A) -- (\O,\M);\draw (\Q,\A) -- (\Q,\M);\draw (\S,\A) -- (\S,\M);
             \draw (\U,\A) -- (\U,\M);\draw (\W,\A) -- (\W,\M);\draw (\Y,\A) -- (\Y,\M);
             \draw (\ZZ,\A) -- (\ZZ,\M);
             \draw (\Za,\A) -- (\Za,\M);
             \draw (\Zb,\A) -- (\Zb,\M);
             \draw (\Zc,\A) -- (\Zc,\M);
             \draw (\Zd,\A) -- (\Zd,\M);
             \draw (\Ze,\A) -- (\Ze,\M);
            \draw[solid,line width=0.4mm] (\A,\A)--(\C,\K) -- (\E,\I) -- (\G,\G) -- (\I,\G) --(\K,\G)-- (\M,\E) -- (\O,\I) -- (\Q,\G) -- (\S,\E) -- (\U,\C)  -- (\W,\A) --(\Y,\E) -- (\ZZ,\E) --(\Za,\I) -- (\Ze, \A);
         \end{tikzpicture}

               \end{center}
         \caption{From left to right, we show a Dyck path $A=UUDDUUUDUUDUDDDDUD$, a Motzkin path  $B=UUFDDFFUDUFUDDFFUD$ and a \L{}ukasiewicz path $C=U_5DDFFDU_2DDDDU_2FU_2DDDD$.}
         \label{fig1}
\end{figure}

Throughout this paper, a {\it lattice path} is defined by a starting point $P_0=(0,0)$, an ending point $P_n=(n,0)$,  consisting of  steps lying in $S=\{(1,i), i\in\mathbb{Z}\}$,
 and never going below the $x$-axis. The {\it length} of a path is the number of its steps. We denote by $\epsilon$ the empty path, {\it i.e.}, the path of length zero.
 Constraining the steps to lie into  $\{(1,1),(1,-1)\}$ (resp. $\{(1,1),(1,0), (1,-1)\}$),   we retrieve the well known definition of Dyck paths (resp. Motzkin paths). {\it \L{}ukasiewicz paths}
 are obtained when the steps belong to $\{(1,i)\in S, i\geq -1\}$. We refer to \cite{Ges,Ran,Sta,Var,Vie} for some combinatorial studies on \L{}ukasiewicz paths. Let $\mathcal{L}_n$, $\mathcal{D}_n$, $\mathcal{M}_n$, $n\geq 0$,
  respectively, be the sets of \L{}ukasiewicz, Dyck and Motzkin paths of length  $n$, and $\mathcal{L}=\cup_{n\geq 0} \mathcal{L}_n$, $\mathcal{D}=\cup_{n\geq 0} \mathcal{D}_n$,
  $\mathcal{M}=\cup_{n\geq 0} \mathcal{M}_n$. For convenience, we set $D=(1,-1)$, $F=(1,0)$, $U=U_1=(1,1)$ and $U_i=(1,i)$ for $i\geq 2$. See Figure \ref{fig1} for an illustration
   of  Dyck, Motzkin and \L{}ukasiewicz paths of length 18. Note that \L{}ukasiewicz paths can be interpreted as an algebraic language  of words $w\in \{x_0,x_1,x_2, \ldots\}^\star$
   such that $\delta(w)=-1$ and $\delta(w')\geq 0$ for any proper prefix $w'$ of $w$ where $\delta$ is the map from $\{x_0,x_1,x_2, \ldots\}^\star$ to $\mathbb{Z}$ defined by $\delta(w_1w_2\ldots w_n)=\sum_{i=1}^n\delta(w_i)$
   with $\delta(x_i)=i-1$ (see \cite{Lot,Sch}).

\begin{figure}[h]
\begin{center}(1)\quad\scalebox{0.55}{\begin{tikzpicture}[ultra thick]
 \draw[black, thick] (0,0)--(9,0); \draw[black, thick] (0,0)--(0,2.5);
  \draw[black, line width=3pt] (0,0)--(0.4,0);
  \draw[orange,very thick] (0.4,0) parabola bend (1.9,2) (3.4,0);
 \draw  (1.8,0.9) node {\huge $L'$};
 \end{tikzpicture}}
 \qquad (2)\quad\scalebox{0.55}{\begin{tikzpicture}[ultra thick]
 \draw[black, thick] (0,0)--(9,0); \draw[black, thick] (0,0)--(0,2.5);
  \draw[black, line width=3pt] (0,0)--(0.4,1.6);
  \draw[black, dashed, very thick] (0.4,1.6)--(1.7,1.6);
   \draw[black, line width=3pt] (1.7,1.6)--(2.1,1.2);
   \draw[black, dashed, very thick] (2.1,1.2)--(3.3,1.2);
   \draw[black, line width=3pt] (3.3,1.2)--(3.7,0.8);
   \draw[black, dashed, very thick] (3.7,0.8)--(4.9,0.8);
    \draw[black, line width=3pt] (4.9,0.8)--(5.3,0.4);
   \draw[black, dashed, very thick] (5.3,0.4)--(6.5,0.4);
   \draw[black, line width=3pt] (6.5,0.4)--(6.9,0);

 \draw[orange,very thick] (0.4,1.6) parabola bend (1.05,2.4) (1.7,1.6);

  \draw[orange,very thick] (2.1,1.2) parabola bend (2.7,2) (3.3,1.2);

  \draw[orange,very thick] (5.3,0.4) parabola bend (5.9,1.2) (6.5,0.4);
  \draw[orange,very thick] (6.9,0) parabola bend (7.7,1.2) (8.5,0);
 \draw  (1.1,1.9) node {$L_1$};\draw  (2.6,1.5) node { $L_2$};
 \draw  (5.85,0.8) node { $L_k$};\draw  (7.7,0.6) node { $L'$};
 \end{tikzpicture}}

\end{center}
\caption{The two forms of the decomposition of  a non-empty \L{}ukasiewicz path.}
\label{fig2}\end{figure}
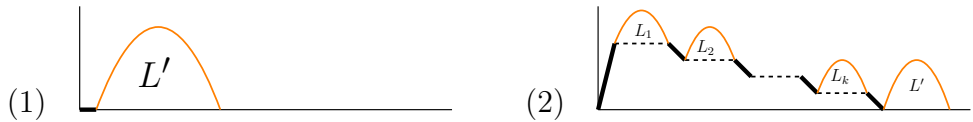

Any non-empty \L{}ukasiewicz path $L\in \mathcal{L}$ can be decomposed (see \cite{Fla}) into one of the two following forms: (1) $L=FL'$ with $L'\in\mathcal{L}$, or (2) $L=U_kL_1DL_2D\ldots L_kD L'$ with $k\geq 1$ and $L_1,L_2, \ldots, L_k, L'\in \mathcal{L}$ (see Figure \ref{fig2}).

 Due to this decomposition, the generating function $L(x)$ for the cardinalities of the sets $\mathcal{L}_n$, $n\geq 0$, satisfies the functional equation $L(x)=1+xL(x)+\sum_{k\geq 1}x^{k+1}L(x)^{k+1}$,
 or equivalently,  $L(x)=\frac{1}{1-xL(x)}$. Then, $L(x)=\frac{1-\sqrt{1-4x}}{2x}$ and the coefficient of $x^n$ in the series expansion of $L(x)$ is given by the $n$-th Catalan
 number $\frac{1}{n+1}{2n \choose n}$ (see sequence \href{https://oeis.org/A000108}{A000108} in \cite{Sloa}).

  A {\it pattern} of length one (resp. two) in a lattice path $L$ consists of one step (resp. two consecutive steps). We will say that an {\it occurrence} of a pattern is at {\it position} $i\geq 1$,  in $L$ whenever the first step of this occurrence appears at the $i$-th step of the path.
   The {\it height} of an occurrence is the minimal ordinate reached by its points.  For instance, the path $C=U_5DDF\mathbf{FD}U_2DDDDU_2FU_2DDDD$ (see Figure \ref{fig1}) contains
   one occurrence of the  pattern $FD$ at position $5$ and of height $2$.

  Following the recent studies \cite{Bar1,Bar2,Man}, we define an equivalence relation on the set $\mathcal{L}$ for a given pattern $\alpha$:
   {\it two \L{}ukasiewicz paths of the same length are $\alpha$-equivalent  whenever the occurrences of the pattern $\alpha$ appear at the same positions in the two paths.}
   For instance, $UFFF\mathbf{FD}UDUDFFFFUDFF$ is $FD$-equivalent to the path $C$ in Figure \ref{fig1}   since the only one occurrence of the pattern $FD$ (in boldface) appear at
   the same position in the two paths. Note that the height of the occurrences of $\alpha$ does not involve in this definition.

 In this paper, for any pattern $\alpha$ of length at most two,  we consider the above equivalence relation on the set $\mathcal{L}$, and for each of them we provide the cardinality of the quotient set with respect to the length.
  Three general methods are used:

    - ($M_0$) we prove that any $\alpha$-equivalence class contains at least one Motzkin path. Using $\mathcal{M}_n\subset \mathcal{L}_n$ for $n\geq 0$, we deduce that the number of $\alpha$-equivalence classes in $\mathcal{M}_n$ is equal to that of $\mathcal{L}_n$. Since  the authors have already determined this number for $\mathcal{M}_n$ (see \cite{Bar2}), we can conclude,

 - ($M_1$) we directly count the number of subsets $\{i_1,i_2,\ldots, i_k\}\subset\{1,2,\ldots, n\}$ of the possible occurrence positions of $\alpha$  in a \L{}ukasiewicz path in $\mathcal{L}_n$,

  - ($M_2$) we exhibit a one-to-one correspondence between a subset of \L{}ukasiewicz paths (subset of representatives of the  classes) and the set of equivalence classes by using combinatorial reasonings, and then, we evaluate algebraically the generating function for this subset.

   The paper is organized as follows. In Section 2, we consider the case of patterns $\alpha$ studied using method ($M_0$), {\it i.e.}, $\alpha\in \{U,UU,UD,UF,DU,FU\}$.
   In Section 3, we focus on these ones that can be dealt using method  ($M_1$), {\it i.e.}, $\alpha\in \{F,D,FD,DF,DD\}$. In Section 4, we complete our study by the remaining cases  which are obtained using method ($M_2$).
   We refer to Table \ref{tab1}  for an exhaustive list of our enumerative results.

\begin{table}[h]\label{tab1}
\begin{center}
\scalebox{0.9}{$\begin{array}{|c|c|c|c|c|}\hline
\mbox{ Pattern } \alpha&\mbox{ Sequence }&\mbox{ Sloane }& a_{n}, 1\leq n\leq 10&\mbox{Method}\\\hline

U& {n\choose\lfloor\frac{n}{2}\rfloor} &\href{https://oeis.org/A001405}{A001405} &1,2,3,6,10,20,35,70,126,252&\multirow{5}*{$M_0$}\\
\cline{1-4}
UU& {\frac {1-2\,x+{x}^{2}-\sqrt { \left( {x}^{2}+1 \right)  \left(
1-3\,{x}^{2} \right) }}{2x(-1+2x-x^2+{x}^{3})}}
 & A191385&1,1,1,2,3,5,7,12,18,31&\\
\cline{1-4}
UD&  \mbox{Fibonacci}&\href{https://oeis.org/A005251}{A005251} &1,2,3,5,8,13,21,34,55,89& \\
\cline{1-4}
UF,FU & \frac{2}{1-2x-\sqrt{1-4x^3}} &\href{https://oeis.org/A165407}{A165407}&1,1,2,3,4,7,11,16,27,43&\\
\cline{1-4}
DU& \mbox{Shift of Fibonacci} &\href{https://oeis.org/A000045}{A000045} &1,1,1,2,3,5,8,13,21,34&\\
\hline\hline
F&2^n-n & \href{https://oeis.org/A000325}{A000325}&1,2,5,12,27,58,121,248,503,1014&\multirow{4}*{$M_1$}\\
\cline{1-4}
D& 2^{n-1} & \href{https://oeis.org/A011782}{A011782} &1,1,2,4,8,16,32,64,128,256&\\
\cline{1-4}
FD,DF&\mbox{Fibonacci} & \href{https://oeis.org/A005251}{A005251} &1,1,2,3,5,8,13,21,34,55&\\
\cline{1-4}
DD&\frac{1-x}{1-2x+x^2-x^3}& \href{https://oeis.org/A005251}{A005251} &1,1,2,4,7,12,21,37,65,114&\\
\hline\hline
U_k& \mbox{Motzkin} & \href{https://oeis.org/A001006}{A001006} &1,2,4,9,21,51,127,323,835,2188&\multirow{5}*{$M_2$}\\
\cline{1-4}
FF& {\frac {1-3x+4x^2-5x^3+7x^4-7x^5+6x^6-3x^7+{x}^{8}}{ \left(1-2x+x^2 -{x}^{3} \right)  \left( 1-
x \right) ^{2}}}& \mbox{New} &1,2,2,5,9,17,32,59,107,192&\\
\cline{1-4}
FU_k,U_kF& \frac{1-x +2x^2+ \sqrt{1-2x+x^2-4x^3}}{1-2x-3x^3+(1-x+x^2)\sqrt{1-2x+x^2-4x^3}} & \href{https://oeis.org/A023431}{A023431} &1,1,2,4,7,13,26,52,104,212&\\
\cline{1-4}
U_kD& \frac{1-x +x^2+ \sqrt{1-2x-x^2-2x^3+x^4}}{1-2x-x^3+(1-x)\sqrt{1-2x-x^2-2x^3+x^4}} & \href{https://oeis.org/A292460}{A292460} &1, 2, 4, 8, 17, 37, 82, 185, 423, 978&\\
\cline{1-4}
DU_k&\frac{1-x -x^2-2x^3+ \sqrt{1-2x-x^2-2x^3+x^4}}{1-2x-x^3+(1-x)\sqrt{1-2x-x^2-2x^3+x^4}}  & \href{https://oeis.org/A004148}{A004148} &1,1,1, 2, 4, 8, 17, 37, 82, 185&\\
\hline
\end{array}$}
\caption{Number of $\alpha$-equivalence classes for \L{}ukasiewicz paths. The last three sequences are recorded in OEIS \cite{Sloa} as generalized Catalan sequences. }
\end{center}
\end{table}

\section{Modulo   $\alpha\in \{U,UU,UD,UF,DU,FU\}$}

In this section, we focus on the patterns that can be dealt using method ($M_0$).
\begin{lem}\label{le1} For $n\geq 0$, let $L$ be a  \L{}ukasiewicz path in $\mathcal{L}_n$ and $\alpha\in\{U,UU,UD,UF,DU,FU\}$. Then, there exists a Motzkin path $M\in\mathcal{M}_n$ such that $M$ and $L$ are $\alpha$-equivalent.\end{lem}
\noindent {\it Proof.} Let us assume that $\alpha\in\{U,UU,UD\}$. Any non-empty \L{}ukasiewicz path can be decomposed into one of the two following forms:  ($i$) $L=FL'$ with $L'\in\mathcal{L}$, and ($ii$) $L=U_kL_1DL_2D\ldots L_kD L'$ with $k\geq 1$ and $L_1,L_2, \ldots, L_k, L'\in \mathcal{L}$.
For $n\geq 0$, we recursively define a map $\phi$ from $\mathcal{L}_n$ to $\mathcal{M}_n$ as follows:

$$\left\{\begin{array}{ll}\phi(\epsilon)&=\epsilon,\\
\phi(FL')&=F\phi(L'),\\
\phi(UL_1DL')&=U\phi(L_1)D\phi(L'),\\
\phi(U_kL_1DL_2D\ldots L_kD L')&=F\phi(L_1)F\phi(L_2)F\ldots\phi(L_k)F\phi(L')\qquad\mbox{  for }k\geq 2.\end{array}\right.$$

 Clearly, $\phi(L)$ is a Motzkin path in $\mathcal{M}_n$, and whenever $\alpha\in\{U,UU,UD\}$ the occurrence positions of $\alpha$ in $L$ and $\phi(L)$ are identical.
  Then, the equivalence class of $L$ contains a Motzkin path $\phi(L)$.

Let us assume that $\alpha=DU$. Any non-empty \L{}ukasiewicz path $L$ can be written as follows: $$L=K_0\prod_{i=1}^r (DU)^{a_i}K_i$$ with $r\geq 0$, $a_i\geq 1$ for $1\leq i\leq r$, and where $K_i$, $0\leq i\leq r$,
 are some parts that do not contain any pattern $DU$. Note that $K_0$ and $K_r$ necessarily contain at least one step. From $L\in\mathcal{L}_n$, we define the Motzkin path
$$M=UF^{b_0-1}\left(\prod_{i=1}^{r-1} (DU)^{a_i} F^{b_i}\right) (DU)^{a_r}DF^{b_r-1}\in\mathcal{M}_n$$ where $b_i=|K_i|$ for $0\leq i\leq r$. Since the occurrence positions of $DU$ in $L$ and $M$ are identical,
$M$ is a Motzkin path in the same class as $L$.

 Let us assume that $\alpha\in\{FU,UF\}$. Any non-empty \L{}ukasiewicz path $L$ can be written as follows: $$L=K_0\prod_{i=1}^r \alpha^{a_i}K_i$$ with $r\geq 0$, $a_i\geq 1$ for $1\leq i\leq r$, and where $K_i$, $0\leq i\leq r$,
 are some parts that  do not contain any pattern $\alpha$. From $L\in\mathcal{L}_n$, we define the Motzkin path
 $$M=F^{b_0}\prod_{i=1}^r \alpha^{a_i} D^{c_i}F^{b_i-c_i}\in\mathcal{M}_n$$ where $b_0=|K_0|$, and for $1\leq i\leq r$, $b_i=|K_i|$ and $c_i=\min\{b_i,a_i+\sum_{j=1}^{i-1}(a_j-c_j)\}$.
 Less formally, $K_0$ is replaced with $F^{|K_0|}$, and for $i$ from $1$ to $r$,  $K_i$ is replaced with  $D^{c_i}F^{|K_i|-c_i}$ where the value $c_i$ is the maximal number of
  down steps $D$ that can be placed so that $M$ remains a lattice path. This ensures that $M$ has the same occurrence positions of $\alpha$ as $L$, which means that $M$ is a Motzkin path in the same class as $L$.\hfill $\Box$

\medskip

Using Lemma \ref{le1} and the fact that $\mathcal{M}\subset\mathcal{L}$, we directly deduce the following theorem.

\begin{thm}\label{thm1} For $\alpha\in\{U,UU,UD,UF,DU,FU\}$ and $n\geq 0$, the number of $\alpha$-equivalence classes in $\mathcal{L}_n$ also is that of  $\mathcal{M}_n$. \end{thm}

Since the authors  have already determined  the number of $\alpha$-equivalence classes in $\mathcal{M}_n$, we refer to their paper \cite{Bar2} for a detailed description of the different proofs,
 and we report the results in Table \ref{tab1}.

\section{Modulo   $\alpha\in \{F,D,FD,DF,DD\}$}

In this section, we focus on the patterns that can be dealt with method ($M_1$) which consists in counting directly the
 possible subsets of occurrence positions of the pattern in  a \L{}ukasiewicz path.

\begin{thm}\label{thm2} The number of $F$-equivalence classes in $\mathcal{L}_n$, $n\geq 0$, is given by $2^n-n$ (see sequence \href{https://oeis.org/A000325}{A000325} in \cite{Sloa}). \end{thm}
\noindent {\it Proof.} Let $L$ be a \L{}ukasiewicz path of length $n\geq 1$, and let $1\leq i_1<i_2< \ldots< i_{\ell}\leq n$, $0\leq \ell\leq n$,  be the sequence of  occurrence positions of $F$ in $L$.  Since a \L{}ukasiewicz path cannot contain exactly $n-1$ occurrences of $F$, we have $\ell\neq n-1$. Now, let us prove that for any $\ell\neq n-1$, $0\leq \ell\leq n$, and for any sequence $1\leq i_1<i_2< \ldots< i_{\ell}\leq n$, there exists $L\in \mathcal{L}_n$ where the positions of its flats are exactly $i_1, \ldots, i_\ell$. We distinguish two cases:  $n-\ell$ is odd 
(different from one), and  $n-\ell$ is even.

  If $n-\ell\neq 1$ is odd, then we define the \L{}ukasiewicz path $K=U_2DD (UD)^k$ of length $n-\ell$ (that is $k=\frac{n-\ell-3}{2}$); otherwise, we define the  \L{}ukasiewicz path $K=(UD)^{\frac{n-\ell}{2}}$.

   For these two cases, we consider the \L{}ukasiewicz  path $L$ (of length $n$) obtained by inserting $\ell$ flats in $K$ so that the positions of flats in $L$  are given by $i_1, \ldots, i_\ell$.

  Then, any increasing sequence $1\leq i_1<i_2<\ldots<i_\ell\leq n$, $0\leq \ell\leq n$ and $n-\ell\neq 1$, corresponds to the positions of $F$ in a \L{}ukasiewicz  path. Since there are $2^n -n$ possible such sequences, the proof is completed.
\hfill $\Box$

\begin{thm}\label{thm3} The number of $D$-equivalence classes in $\mathcal{L}_n$, $n\geq 0$,  is given by $2^{n-1}$ (see sequence \href{https://oeis.org/A011782}{A011782} in \cite{Sloa}). \end{thm}
\noindent {\it Proof.} Let $L$ be a \L{}ukasiewicz path of length $n\geq 1$, and let $1\leq i_1<i_2< \ldots< i_{\ell}\leq n$,
$0\leq \ell\leq n$,  be the sequence of occurrence positions  of $D$ in $L$.  Since the first step of a \L{}ukasiewicz
 path cannot be a down step $D$,  we have $0\leq\ell\leq n-1$ and $i_1\neq 1$. Now, let us prove that for any $\ell$, $0\leq\ell\leq n-1$,
  and for any sequence $2\leq i_1<i_2< \ldots< i_{\ell}\leq n$, there exists $L\in \mathcal{L}_n$ where the positions of its down steps
   are exactly $i_1, \ldots, i_\ell$.

We define the \L{}ukasiewicz path $L$ as follows:
$$L= U_{\ell}F^{i_1-2} \prod_{j=1}^\ell D F^{i_{j+1}-i_j-1},$$
where $i_{\ell+1}=n+1$.

This definition  ensures that the down steps appear on positions $i_1,i_2, \ldots , i_\ell$. Therefore, any increasing sequence $2\leq i_1<i_2<\ldots<i_\ell\leq n$,
$0\leq \ell\leq n-1$, corresponds to the sequence of positions of $D$ in a \L{}ukasiewicz  path. Since there are $2^{n-1}$ possible such sequences, the proof is completed.\hfill $\Box$

\begin{thm}\label{thm4} The number of $FD$-equivalence (resp. $DF$-equivalence) classes in $\mathcal{L}_n$, $n\geq 0$, is given by the Fibonacci number defined by
$f_0=1$, $f_1=1$, $f_2=1$ and $f_n=f_{n-1}+f_{n-2}$ for $n\geq 3$ (see sequence \href{https://oeis.org/A005251}{A005251} in \cite{Sloa}). \end{thm}
\noindent {\it Proof.}  Let $L$ be a  \L{}ukasiewicz path of length $n\geq 1$. For $\alpha\in \{FD,DF\}$,
let $1\leq i_1<i_2< \ldots< i_{\ell}\leq n$, $0\leq \ell\leq n$,  be the sequence of  occurrence positions  of $\alpha$ in $L$.
 Since a \L{}ukasiewicz path cannot contain two occurrences of $\alpha$ at two adjacent positions $i$ and $i+1$, we necessarily
  have $i_j-i_{j-1}>1$ for $j\geq 2$. A pattern $\alpha$ cannot appear at position one, which implies $i_1\geq 2$. Also, a pattern $\alpha$ cannot appear at position $n$, which implies $i_{\ell}\leq n-1$.
    Now, let us prove that for any $\ell$, $0\leq \ell\leq n$, and for any sequence $2\leq i_1<i_2< \ldots< i_{\ell}\leq n-1$ with
   $i_j-i_{j-1}>1$ for $j\geq 2$, there exists $L\in \mathcal{L}_n$ where the positions of its $\alpha$ are exactly $i_1, \ldots, i_\ell$.
   We set $I=\{i_1, \ldots, i_\ell\}$. So we define the \L{}ukasiewicz path $L$ as follows:
  $$L= U_{\ell}F^{i_1-2} \prod_{j=1}^\ell \alpha F^{i_{j+1}-i_j-2},$$
where $i_{\ell+1}=n+1$.

This definition  ensures that the down steps appear on positions $i_1,i_2, \ldots , i_\ell$.
  Therefore, any increasing sequence $2\leq i_1<i_2<\ldots<i_\ell\leq n-1$, $0\leq \ell\leq n$ with $i_j-i_{j-1}>1$ for $j\geq 2$, corresponds to the sequence
   of positions of $\alpha$ in a \L{}ukasiewicz  path. It is well known (see \cite{Vajd}  for instance) that such a sequence is enumerated by the Fibonacci number $f_n$
   defined by $f_1=1$, $f_2=1$, and $f_n=f_{n-1}+f_{n-2}$ for $n\geq 3$.\hfill $\Box$

\begin{thm}\label{thm5} The number of $DD$-equivalence  classes in $\mathcal{L}_n$, $n\geq 0$, is given by the $n$-th term $g_n$ of the sequence  defined by  $g_0=1$, $g_1=1$, $g_2=1$, $g_3=2$ and $g_n=g_{n-1}+g_{n-2}+g_{n-4}$ for $n\geq 4$ (see sequence
 \href{https://oeis.org/A005251}{A005251}  in \cite{Sloa}). \end{thm}
\noindent {\it Proof.}  Let $L$ be a  \L{}ukasiewicz path of length $n\geq 1$, and let $1\leq i_1<i_2< \ldots< i_{\ell}\leq n$,
 $0\leq \ell\leq n$,  be the sequence of  occurrence positions  of $DD$ in $L$.  Whenever a \L{}ukasiewicz path contains two
  occurrences of $DD$ at two positions $i$ and $i+2$, $L$ necessarily have an occurrence of $DD$ at position $i+1$,
  that is $i_{j+1}\neq i_j+2$. A pattern $DD$ cannot appear at position one, which implies $i_1\geq 2$. Also, a pattern $DD$ cannot appear at position $n$, which implies $i_\ell\leq n-1$.
  Now, let us prove that
   for any $\ell$, $0\leq \ell\leq n$, and for any sequence $2\leq i_1<i_2< \ldots< i_{\ell}\leq n-1$ satisfying $i_j-i_{j-1}\neq 2$ for
    $j\geq 2$, there exists $L\in \mathcal{L}_n$ where the positions of its $DD$ are exactly $i_1, \ldots, i_\ell$.
     We set $I=\{i_1, \ldots, i_\ell\}$, $I^{+}=\{i_1+1, \ldots, i_\ell+1\}$. We consider  the unique partition $G_1, G_2, \ldots , G_r$, $r\geq 1$,
      of $[1,n+1]\backslash (I\cup I^{+})$ such that $G_i$, $1\leq i\leq r$, is a maximal non empty interval satisfying $\max{G_i}<\min{G_{i+1}}$ for $i\leq r-1$.
      So we define the \L{}ukasiewicz path $L$ as follows:
      $$L= U_{b}F^{i_1-2} D^{a_1} \left(\prod_{j=2}^{r-1} F^{g_j}D^{a_j}\right)   F^{|G_{r}|-1}     $$
where $b=|I\cup I^{+}|-1$ is the number of down steps $D$ in $L$,  $a_j=\min{G_{j+1}}-\max{G_{j}}-1$  and $g_j=|G_j|$ for $1\leq j\leq r-1$. Less formally, we place occurrences $DD$ on positions $i_1,\ldots , i_\ell$, we start the path with $U_b$ where $b$ is the
number of down steps $D$, and we place flat steps anywhere else.

For instance, if $n=14$ and $I=\{2,3,7,10,13\}$, then we have $[1,15]\backslash (I\cup I^{+})=\{1,5,6,9,12,15\}$ $b=15-6=9$ and $G_1=\{1\}$, $G_2=\{5,6\}$, $G_3=\{9\}$, $G_4=\{12\}$, $G_5=\{15\}$,
 which induces $L=U_9 DDD FF DD F DD F DD$.

 Therefore, any increasing sequence $2\leq i_1<i_2<\ldots<i_\ell\leq n-1$, $0\leq \ell\leq n$ with $i_j-i_{j-1}\neq 2$ for $j\geq 2$, corresponds to the set of positions
 of $DD$ in a \L{}ukasiewicz  path. It is already known (see \cite{Aus} for instance) that such sequences are enumerated by the general term $g_n$ of \href{https://oeis.org/A005251}{A005251} in \cite{Sloa} defined by
  $g_0=1$, $g_1=1$, $g_2=1$, $g_3=2$ and $g_n=g_{n-1}+g_{n-2}+g_{n-4}$ for $n\geq 4$.\hfill $\Box$

\section{Other patterns}

 In this section, we consider the  equivalence relation  on $\mathcal{L}$ where two paths $L$ and $L'$ belong to the same class whenever  for any $k\geq 1$,
 the  occurrence positions  of $U_k$ (resp.  $DU_k$, $U_kD$, $FU_k$, $U_kF$) are the same in $L$ and $L'$. Also, we study the $FF$-equivalence relation in $\mathcal{L}$.
 For all these cases, we use the method ($M_2$) that consists in exhibiting subsets of representatives of equivalence classes, and
 determining algebraically their cardinalities.

\subsection{Modulo the up steps $U_k$, $k\geq 1$}
     In this part, we assume that two \L{}ukasiewicz paths $L$ and $L'$ of the same length are $U_k$-equivalent whenever for any $k\geq 1$, $L$ and $L'$  have the same  positions of $U_k$.

Let $\mathcal{B}$ be the set of \L{}ukasiewicz paths without any flat steps at positive height. For instance, we have $U_3DDDFUD\in\mathcal{B}$ and $U_3FDDDUD\notin\mathcal{B}$.
Let $\bar{\mathcal{B}}\subset \mathcal{B}$ be the set of \L{}ukasiewicz paths without any flat steps.

\begin{lem}\label{le2}
There is a bijection between $\mathcal{B}$ and the set of $U_k$-equivalence classes of $\mathcal{L}$.
\end{lem}

\begin{figure}[h]

 \begin{center} $L=$\begin{tabular}{c}\begin{tikzpicture}[scale=0.1]
            \draw (\A,\A)-- (38,\A);
             \draw (\A,\E)-- (\Ze,\E);
              \draw (\A,\C)-- (\Ze,\C);
               \draw (\A,\G)-- (\Ze,\G);
               \draw (\A,\I)-- (\Ze,\I);
              \draw (\A,\K)-- (\Ze,\K);
             \draw (\A,\M)-- (\Ze,\M);
            \draw (\A,\A) -- (\A,\O);
             \draw (\C,\A) -- (\C,\M);\draw (\E,\A) -- (\E,\M);\draw (\G,\A) -- (\G,\M);
             \draw (\I,\A) -- (\I,\M);\draw (\K,\A) -- (\K,\M);\draw (\M,\A) -- (\M,\M);
             \draw (\O,\A) -- (\O,\M);\draw (\Q,\A) -- (\Q,\M);\draw (\S,\A) -- (\S,\M);
             \draw (\U,\A) -- (\U,\M);\draw (\W,\A) -- (\W,\M);\draw (\Y,\A) -- (\Y,\M);
             \draw (\ZZ,\A) -- (\ZZ,\M);
             \draw (\Za,\A) -- (\Za,\M);
             \draw (\Zb,\A) -- (\Zb,\M);
             \draw (\Zc,\A) -- (\Zc,\M);
             \draw (\Zd,\A) -- (\Zd,\M);
             \draw (\Ze,\A) -- (\Ze,\M);
            \draw[solid,color=red, line width=0.4mm] (\A,\A)--(\C,\G);
             \draw[solid, line width=0.4mm] (\C,\G) -- (\E,\E);
              \draw[solid, color=red,line width=0.4mm] (\E,\E)-- (\G,\G) ;
               \draw[solid, line width=0.4mm] (\G,\G) -- (\I,\E);
            \draw[solid,line width=0.4mm] (\I,\E)  -- (\K,\E);
             \draw[solid,line width=0.4mm] (\K,\E) -- (\M,\E) -- (\O,\E) -- (\Q,\E);
            \draw[solid,color=red, line width=0.4mm] (\Q,\E)--(\S,\G) -- (\U,\I);
             \draw[solid,line width=0.4mm] (\U,\I) -- (\W,\G) --(\Y,\E) -- (\ZZ,\E) --(\Za,\E);
             \draw[solid,line width=0.4mm] (\Za,\E) -- (\Zb,\C);
            \draw[solid,line width=0.4mm] (\Zb,\C) -- (\Zc,\A) -- (\Zd,\A) -- (\Ze, \A);
         \end{tikzpicture}\\
         \end{tabular}
$\longrightarrow L'= $
          \begin{tabular}{c}\begin{tikzpicture}[scale=0.1]
            \draw (\A,\A)-- (38,\A);
             \draw (\A,\E)-- (\Ze,\E);
              \draw (\A,\C)-- (\Ze,\C);
               \draw (\A,\G)-- (\Ze,\G);
               \draw (\A,\I)-- (\Ze,\I);
              \draw (\A,\K)-- (\Ze,\K);
               \draw (\A,\M)-- (\Ze,\M);
            \draw (\A,\A) -- (\A,\O);
             \draw (\C,\A) -- (\C,\M);\draw (\E,\A) -- (\E,\M);\draw (\G,\A) -- (\G,\M);
             \draw (\I,\A) -- (\I,\M);\draw (\K,\A) -- (\K,\M);\draw (\M,\A) -- (\M,\M);
             \draw (\O,\A) -- (\O,\M);\draw (\Q,\A) -- (\Q,\M);\draw (\S,\A) -- (\S,\M);
             \draw (\U,\A) -- (\U,\M);\draw (\W,\A) -- (\W,\M);\draw (\Y,\A) -- (\Y,\M);
             \draw (\ZZ,\A) -- (\ZZ,\M);
             \draw (\Za,\A) -- (\Za,\M);
             \draw (\Zb,\A) -- (\Zb,\M);
             \draw (\Zc,\A) -- (\Zc,\M);
             \draw (\Zd,\A) -- (\Zd,\M);
             \draw (\Ze,\A) -- (\Ze,\M);
 \draw[solid,color=red, line width=0.4mm] (\A,\A)--(\C,\G);
             \draw[solid, line width=0.4mm] (\C,\G) -- (\E,\E);
              \draw[solid, color=red,line width=0.4mm] (\E,\E)-- (\G,\G) ;
            \draw[solid,line width=0.4mm]  (\G,\G) -- (\I,\E) -- (\K,\C);
            \draw[solid,line width=0.4mm] (\K,\C)-- (\M,\A);
             \draw[solid,line width=0.4mm] (\M,\A)-- (\O,\A) -- (\Q,\A);
            \draw[solid,color=red, line width=0.4mm]  (\Q,\A) -- (\S,\C) -- (\U,\E);
             \draw[solid,line width=0.4mm] (\U,\E) -- (\W,\C) --(\Y,\A) -- (\ZZ,\A) --(\Za,\A) -- (\Zb,\A);
           \draw[solid,line width=0.4mm]  (\Zb,\A)-- (\Zc,\A);
             \draw[solid,line width=0.4mm] (\Zc,\A)-- (\Zd,\A) -- (\Ze, \A);
         \end{tikzpicture}
         \\
         \end{tabular}\\
\end{center}

         \caption{Illustration of the  example described in the proof of Lemma \ref{le2}.}
         \label{fig3}
\end{figure}
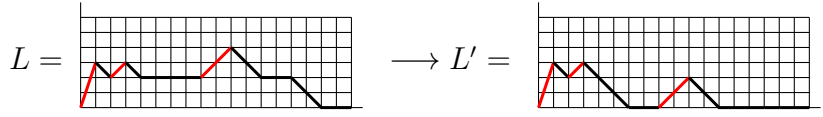

\noindent {\it Proof.} Let $L$ be a non-empty \L{}ukasiewicz path in $\mathcal{L}$. Let us prove that there exists a \L{}ukasiewicz path $L'\in\mathcal{B}$ (with the same length as $L$)
such that $L$ and $L'$ are equivalent. We write $$L=K_0\prod_{i=1}^r \alpha_iK_i$$ with $r\geq 0$,
  where $K_i$ is a part that does not contain any up steps for $0\leq i\leq r$, and $\alpha_i\in\{U_k, k\geq 1\}$ for $1\leq i\leq r$. From $L\in\mathcal{L}$, we define the
   \L{}ukasiewicz path
 $$L'=F^{b_0}\prod_{i=1}^r \alpha_i D^{c_i}F^{b_i-c_i}$$ with $b_0=|K_0|$, and for $1\leq i\leq r$, $b_i=|K_i|$, $c_i=\min\{b_i,a_i+\sum_{j=1}^{i-1}(a_j-c_j)\}$  where $\alpha_i=U_{a_i}$.
  Less formally, $K_0$ is replaced with $F^{|K_0|}$ and for $i$ from $1$ to $r$, $K_i$ is replaced with $D^{c_i}F^{|K_i|-c_i}$ where $c_i$ is the maximal
   number of down steps that can be placed  so that $L'$ remains a  \L{}ukasiewicz path.
 Clearly, $L'$ belongs to $\mathcal{B}$ (it does not contain any flat at positive height), and for any $k\geq 1$ the  occurrence positions
  of $U_k$ are the same as for $L$, {\it i.e.}, $L'\in\mathcal{B}$ is in the same class as $L$.
   For instance, if $L=U_3DUDFFFFUUDDFFDDFF$, then we obtain $L'=U_3DUDDDFFUUDDFFFFFF$ (see Figure \ref{fig3} for an illustration of this example).

      Since the positions of the up steps $U_k$, $k\geq 1$, remain fixed inside a class, and that any flat of $L'\in \mathcal{B}$ lies  necessarily on the $x$-axis,
      there are no other paths in $\mathcal{B}$ in the same class as $L$. The proof is completed.\hfill $\Box$

\begin{thm}\label{thm6} The generating function for the set of  $U_k$-equivalence classes of $\mathcal{L}$ with respect to the length is given by $$\frac{1-x - \sqrt{1-2x-3x^2}}{2x^2},$$ which generates the Motzkin numbers
(\href{https://oeis.org/A001006}{A001006} in \cite{Sloa}). \end{thm}
\noindent {\it Proof.} Using Lemma \ref{le2}, it suffices to obtain the generating function $B(x)$ for the set $\mathcal{B}$. A non-empty \L{}ukasiewicz path $L\in\mathcal{B}$ can be written either $L=FL'$ where $L'\in \mathcal{B}$,
 or $L=U_kL_1 D L_2 D \ldots L_kDL'$ for $k\geq 1$ where
and $L_1, L_2, \ldots, L_k\in\bar{\mathcal{B}}$ are some  \L{}ukasiewicz paths without flats, and $L'\in \mathcal{B}$. So we obtain the functional equation
$B(x)=1+xB(x)+xB(x)\sum_{k\geq 0}x^k\bar{B}(x)^k$ where $\bar{B}(x)$ is the generating function for the set $\bar{\mathcal{B}}$ of \L{}ukasiewicz paths without flats. Using the classical
decomposition of a \L{}ukasiewicz path, $\bar{B}(x)$ satisfies  $\bar{B}(x)=1+\sum_{k\geq 2}x^k\bar{B}(x)^k$, or equivalently $\bar{B}(x)=\frac{1}{(1+x)(1-x\bar{B}(x))}$. A simple calculation provides the result.
\hfill $\Box$
\bigskip

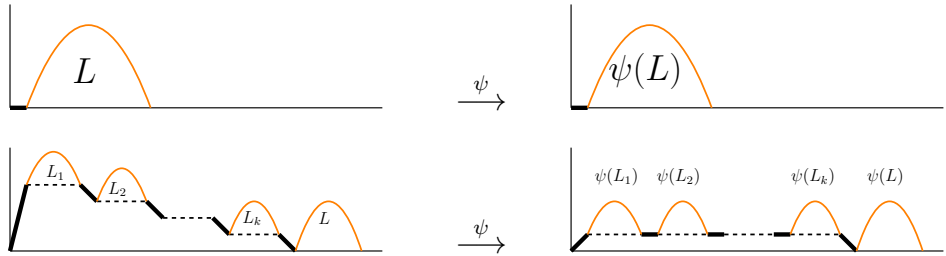
\begin{figure}[h]
\begin{center}\scalebox{0.55}{\begin{tikzpicture}[ultra thick]
 \draw[black, thick] (0,0)--(9,0); \draw[black, thick] (0,0)--(0,2.5);
  \draw[black, line width=3pt] (0,0)--(0.4,0);
  \draw[orange,very thick] (0.4,0) parabola bend (1.9,2) (3.4,0);
 
 \draw  (1.8,0.9) node {\huge $L$};
 \end{tikzpicture}}
 \qquad$\stackrel{\psi}{\longrightarrow}$\qquad
\scalebox{0.55}{\begin{tikzpicture}[ultra thick]
 \draw[black, thick] (0,0)--(9,0); \draw[black, thick] (0,0)--(0,2.5);
  \draw[black, line width=3pt] (0,0)--(0.4,0);
  \draw[orange,very thick] (0.4,0) parabola bend (1.9,2) (3.4,0);
 
 \draw  (1.8,0.9) node {\huge $\psi(L)$};
 \end{tikzpicture}}
\end{center}
\begin{center}\scalebox{0.55}{\begin{tikzpicture}[ultra thick]
 \draw[black, thick] (0,0)--(9,0); \draw[black, thick] (0,0)--(0,2.5);
  \draw[black, line width=3pt] (0,0)--(0.4,1.6);
  \draw[black, dashed, very thick] (0.4,1.6)--(1.7,1.6);
   \draw[black, line width=3pt] (1.7,1.6)--(2.1,1.2);
   \draw[black, dashed, very thick] (2.1,1.2)--(3.3,1.2);
   \draw[black, line width=3pt] (3.3,1.2)--(3.7,0.8);
   \draw[black, dashed, very thick] (3.7,0.8)--(4.9,0.8);
    \draw[black, line width=3pt] (4.9,0.8)--(5.3,0.4);
   \draw[black, dashed, very thick] (5.3,0.4)--(6.5,0.4);
   \draw[black, line width=3pt] (6.5,0.4)--(6.9,0);

 \draw[orange,very thick] (0.4,1.6) parabola bend (1.05,2.4) (1.7,1.6);

  \draw[orange,very thick] (2.1,1.2) parabola bend (2.7,2) (3.3,1.2);

  \draw[orange,very thick] (5.3,0.4) parabola bend (5.9,1.2) (6.5,0.4);
  \draw[orange,very thick] (6.9,0) parabola bend (7.7,1.2) (8.5,0);
 
 \draw  (1.1,1.9) node {$L_1$};\draw  (2.6,1.5) node { $L_2$};
 \draw  (5.85,0.8) node { $L_k$};\draw  (7.6,0.8) node { $L$};
 \end{tikzpicture}}
 \qquad $\stackrel{\psi}{\longrightarrow}$\qquad
\scalebox{0.55}{\begin{tikzpicture}[ultra thick]
 \draw[black, thick] (0,0)--(9,0); \draw[black, thick] (0,0)--(0,2.5);
  \draw[black, line width=3pt] (0,0)--(0.4,0.4);
  \draw[black, dashed, very thick] (0.4,0.4)--(1.7,0.4);
   \draw[black, line width=3pt] (1.7,0.4)--(2.1,0.4);
   \draw[black, dashed, very thick] (2.1,0.4)--(3.3,0.4);
   \draw[black, line width=3pt] (3.3,0.4)--(3.7,0.4);
   \draw[black, dashed, very thick] (3.7,0.4)--(4.9,0.4);
    \draw[black, line width=3pt] (4.9,0.4)--(5.3,0.4);
   \draw[black, dashed, very thick] (5.3,0.4)--(6.5,0.4);
   \draw[black, line width=3pt] (6.5,0.4)--(6.9,0);

 \draw[orange,very thick] (0.4,0.4) parabola bend (1.05,1.2) (1.7,0.4);

  \draw[orange,very thick] (2.1,0.4) parabola bend (2.7,1.2) (3.3,0.4);

  \draw[orange,very thick] (5.3,0.4) parabola bend (5.9,1.2) (6.5,0.4);
  \draw[orange,very thick] (6.9,0) parabola bend (7.7,1.2) (8.5,0);
 
 \draw  (1.1,1.8) node {$\psi(L_1)$};\draw  (2.6,1.8) node {$\psi(L_2)$};
 \draw  (5.85,1.8) node { $\psi(L_k)$};\draw  (7.6,1.8) node { $\psi(L)$};
 \end{tikzpicture}}
\end{center}
\caption{Illustration of the map $\psi$  from $\mathcal{L}$ to $\mathcal{M}$.}
\label{fig4}\end{figure}

Let us define recursively a map $\psi$ from $\mathcal{L}$ to the set of Motzkin paths  $\mathcal{M}$ as follows:
$$\left\{\begin{array}{ll}
\psi(\epsilon)&=\epsilon,\\
\psi(FL)&=F\psi(L),\\
\psi(U_kL_1DL_2D\ldots L_kDL)&=  U\psi(L_1)F\psi(L_2)F\ldots \psi(L_k)D\psi(L),
\end{array}\right.$$ where $L, L_1,L_2, \ldots, L_k$ are are some \L{}ukasiewicz paths.
See Figure \ref{fig4} for an illustration of the bijection $\psi$. For instance, the image by $\psi$ of $U_4FU_2DFDDDU_2U_1DDDFDDFU_2FDU_2DDD$ is
$UFUFFDFFUUDFDFFDFUFFUFDD$. Obviously,  the map $\psi$ preserves the length of the paths.

We easily deduce the two following facts.

\begin{fac}  If $L,L'\in \mathcal{L}$, then $LL'\in \mathcal{L}$ and  we have $\psi(LL')=\psi(L)\psi(L').$
\label{fac1}
\end{fac}
\begin{fac}  There is a one-to-one correspondence between:

\begin{itemize}
\item[(a)]  steps $\{U_k, k\geq 1\}$ in $L$ and   steps $U$ in $\psi(L)$;

\item[(b)] $\{U_kU_\ell, k,\ell\geq 1\}$ in $L$ and $UU$ in $\psi(L)$;

\item[(c)] steps $F$ on the $x$-axis in $L$ and   steps $F$ on the $x$-axis in $\psi(L)$;

\item[(d)] $\{U_kD, k\geq 2\}\cup\{U_kF, k\geq 1\}$  in $L$ and $UF$ in $\psi(L)$;

\item[(e)] peaks $UD$ in $L$ and  peaks $UD$ in $\psi(L)$.

\end{itemize}
\label{fac2}
\end{fac}

\begin{thm}\label{thm6bis} For any $n\geq 0$, the map $\psi$ induces a bijection from $\mathcal{B}_n$ to $\mathcal{M}_n$. \end{thm}
\noindent {\it Proof.} We proceed by induction on $n$. Obviously, for  $n=0$ we have $\psi(\epsilon)=\epsilon$. We assume that $\psi$ is a
bijection from $\mathcal{B}_k$ to $\mathcal{M}_k$, $0\leq k\leq n$, and we prove the result for $n+1$. Using the enumerating result of Theorem \ref{thm6},
 it suffices to prove that $\psi$ is surjective. So, let $M$ be a Motzkin path in  $\mathcal{M}_{n+1}$. We distinguish two cases: ($i$) $M=FM'$ with $M'\in \mathcal{M}_n$, and
 ($ii$) $M=UM'DM''$  where $M'$ and $M''$ are two Motzkin paths in $\mathcal{M}$.

 ($i$) Using the recurrence hypothesis, there is $L'\in \mathcal{B}_n$ such that $M'=\psi(L')$. So, the \L{}ukasiewicz path
 $L=FL'$ lies into $\mathcal{B}_{n+1}$ and  satisfies $\psi(L)=M$ which proves that $M$ belongs to the image by $\psi$ of $\mathcal{B}_{n+1}$.

 ($ii$) We suppose $M=UM'DM''$. We can uniquely write $M'=M_0\prod_{i=1}^rFM_i$ with $r\geq 0$ and where $M_i$ is a (possibly empty) Motzkin path without flat $F$
 on the $x$-axis. Using the recurrence hypothesis, there are $B_0, B_1, \ldots, B_r\in \mathcal{B}$ such that $\psi(B_i)=M_i$, $0\leq i\leq r$.
  Also let $B\in \mathcal{B}$ such that $\psi(B)=M''$.
 Since $B_i$ (resp. $B$) belongs to $\mathcal{B}$, it does not contain any flat at positive height. Since $M_i=\psi(B_i)$, Fact 2(c) implies that
  $B_i$  does not contain any flat on the $x$-axis. So, $B_i$ does not contain any flat steps. So, let us define $$L=U_{r+1}B_0D \left(\prod_{i=1}^rB_iD \right) B.$$
   Clearly, $L$ lies
   in $\mathcal{B}_{n+1}$ and satisfies $\psi(L)=M$; then,  $M$ belongs to the image by $\psi$ of
 $\mathcal{B}_{n+1}$.

The map $\psi$ from $\mathcal{B}_n$ to $\mathcal{M}_n$ is a bijection.\hfill $\Box$

\subsection{Modulo $U_kD$ for $k\geq 1$, and  $U_kF$ for $k\geq 1$  }
      For a given $\alpha\in \{D,F\}$, we define the $U_k\alpha$-equivalence in $\mathcal{L}$ as follows: two \L{}ukasiewicz paths $L$ and $L'$ of the same length are $U_k\alpha$-equivalent whenever for any $k\geq 1$, $L$ and $L'$  have the same  positions of $U_k\alpha$.

Let $\mathcal{C}\subset\mathcal{B}$ be the set of \L{}ukasiewicz paths without any flat steps at positive height and such that any up step $U_k$, $k\geq 1$,
 is immediately followed by a down step $D$. For instance, we have $U_3DDDFUD\in\mathcal{C}$ and $U_3FDDDUD\notin\mathcal{C}$.
 Let $\bar{\mathcal{C}}\subset\mathcal{C}$ be the set of \L{}ukasiewicz paths without flats in $\mathcal{C}$.

\begin{lem}\label{le3}
There is a bijection between $\mathcal{C}$ and the set of $U_kD$-equivalence classes of $\mathcal{L}$.
\end{lem}

\noindent {\it Proof.}
Let $L$ be a non-empty \L{}ukasiewicz path in $\mathcal{L}$. Let us prove that there exists a \L{}ukasiewicz path $L'\in\mathcal{C}$
(with the same length as $L$) such that $L$ and $L'$ belong to the same class. We write
 $$L=K_0\prod_{i=1}^r\left(U_{k_i}DK_i\right),$$ where $k_i\geq 1$ for $1\leq i \leq r$, and $K_0,K_1,K_2, \ldots, K_r$, $r\geq 0$,  are some parts (possibly empty)
 without pattern $U_kD$ for any $k\geq 1$.

We define the \L{}ukasiewicz path $$L'=F^{b_0}\prod_{i=1}^r\left(U_{k_1}DD^{a_i}F^{b_i-a_i}\right),$$ with $b_i=|K_i|$, $0\leq i\leq r$, and  for $1\leq i\leq r$,
 $a_i=\min\{b_i, k_i-1+\sum_{j=1}^{i-1} (k_j-1-a_j)\}$.
Less formally, $K_0$ is replaced with $F^{|K_0|}$, and $K_i$ is replaced with $D^{a_i}F^{b_i-a_i}$ where the value $a_i$, $1\leq i\leq r$, is the maximal number of down steps that can be placed between
 the two occurrences $U_{k_{i}}D$ and $U_{k_{i+1}}D$ so that $L'$ remains a lattice path. Clearly, $L'\in\mathcal{C}$ and $L'$ belongs to the same class as $L$.

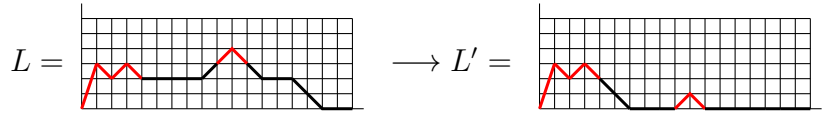
\begin{figure}[h]

 \begin{center} $L=$\begin{tabular}{c}\begin{tikzpicture}[scale=0.1]
            \draw (\A,\A)-- (38,\A);
             \draw (\A,\E)-- (\Ze,\E);
              \draw (\A,\C)-- (\Ze,\C);
               \draw (\A,\G)-- (\Ze,\G);
               \draw (\A,\I)-- (\Ze,\I);
              \draw (\A,\K)-- (\Ze,\K);
             \draw (\A,\M)-- (\Ze,\M);
            \draw (\A,\A) -- (\A,\O);
             \draw (\C,\A) -- (\C,\M);\draw (\E,\A) -- (\E,\M);\draw (\G,\A) -- (\G,\M);
             \draw (\I,\A) -- (\I,\M);\draw (\K,\A) -- (\K,\M);\draw (\M,\A) -- (\M,\M);
             \draw (\O,\A) -- (\O,\M);\draw (\Q,\A) -- (\Q,\M);\draw (\S,\A) -- (\S,\M);
             \draw (\U,\A) -- (\U,\M);\draw (\W,\A) -- (\W,\M);\draw (\Y,\A) -- (\Y,\M);
             \draw (\ZZ,\A) -- (\ZZ,\M);
             \draw (\Za,\A) -- (\Za,\M);
             \draw (\Zb,\A) -- (\Zb,\M);
             \draw (\Zc,\A) -- (\Zc,\M);
             \draw (\Zd,\A) -- (\Zd,\M);
             \draw (\Ze,\A) -- (\Ze,\M);
            \draw[solid,color=red, line width=0.4mm] (\A,\A)--(\C,\G);
             \draw[solid, color=red,line width=0.4mm] (\C,\G) -- (\E,\E);
              \draw[solid, color=red,line width=0.4mm] (\E,\E)-- (\G,\G) ;
               \draw[solid, color=red,line width=0.4mm] (\G,\G) -- (\I,\E);
            \draw[solid,line width=0.4mm] (\I,\E)  -- (\K,\E);
             \draw[solid,line width=0.4mm] (\K,\E) -- (\M,\E) -- (\O,\E) -- (\Q,\E);
            \draw[solid, line width=0.4mm] (\Q,\E)--(\S,\G) ;
             \draw[solid,color=red,line width=0.4mm] (\S,\G) --(\U,\I) -- (\W,\G);
             \draw[solid,line width=0.4mm]  (\W,\G)--(\Y,\E) -- (\ZZ,\E) --(\Za,\E);
             \draw[solid,line width=0.4mm] (\Za,\E) -- (\Zb,\C);
            \draw[solid,line width=0.4mm] (\Zb,\C) -- (\Zc,\A) -- (\Zd,\A) -- (\Ze, \A);
         \end{tikzpicture}\\
         \end{tabular}
$\longrightarrow L'= $
          \begin{tabular}{c}\begin{tikzpicture}[scale=0.1]
            \draw (\A,\A)-- (38,\A);
             \draw (\A,\E)-- (\Ze,\E);
              \draw (\A,\C)-- (\Ze,\C);
               \draw (\A,\G)-- (\Ze,\G);
               \draw (\A,\I)-- (\Ze,\I);
              \draw (\A,\K)-- (\Ze,\K);
               \draw (\A,\M)-- (\Ze,\M);
            \draw (\A,\A) -- (\A,\O);
             \draw (\C,\A) -- (\C,\M);\draw (\E,\A) -- (\E,\M);\draw (\G,\A) -- (\G,\M);
             \draw (\I,\A) -- (\I,\M);\draw (\K,\A) -- (\K,\M);\draw (\M,\A) -- (\M,\M);
             \draw (\O,\A) -- (\O,\M);\draw (\Q,\A) -- (\Q,\M);\draw (\S,\A) -- (\S,\M);
             \draw (\U,\A) -- (\U,\M);\draw (\W,\A) -- (\W,\M);\draw (\Y,\A) -- (\Y,\M);
             \draw (\ZZ,\A) -- (\ZZ,\M);
             \draw (\Za,\A) -- (\Za,\M);
             \draw (\Zb,\A) -- (\Zb,\M);
             \draw (\Zc,\A) -- (\Zc,\M);
             \draw (\Zd,\A) -- (\Zd,\M);
             \draw (\Ze,\A) -- (\Ze,\M);
\draw[solid,color=red, line width=0.4mm] (\A,\A)--(\C,\G);
             \draw[solid, color=red,line width=0.4mm] (\C,\G) -- (\E,\E);
              \draw[solid, color=red,line width=0.4mm] (\E,\E)-- (\G,\G) ;
               \draw[solid, color=red,line width=0.4mm] (\G,\G) -- (\I,\E);
            \draw[solid,line width=0.4mm]  (\I,\E) -- (\K,\C);
            \draw[solid,line width=0.4mm] (\K,\C)-- (\M,\A);
             \draw[solid,line width=0.4mm] (\M,\A)-- (\O,\A) -- (\Q,\A);
             \draw[solid, line width=0.4mm]  (\Q,\A) -- (\S,\A);
            \draw[solid,color=red, line width=0.4mm] (\S,\A)  --(\U,\C)--(\W,\A);
             \draw[solid,line width=0.4mm]  (\W,\A) --(\Y,\A) -- (\ZZ,\A) --(\Za,\A) -- (\Zb,\A);
           \draw[solid,line width=0.4mm]  (\Zb,\A)-- (\Zc,\A);
             \draw[solid,line width=0.4mm] (\Zc,\A)-- (\Zd,\A) -- (\Ze, \A);
         \end{tikzpicture}
         \\
         \end{tabular}\\
\end{center}

         \caption{Illustration of the  example described in the proof of Lemma \ref{le3}.}
         \label{fig3bis}
\end{figure}

For instance, from  $L=U_3DUDFFFFUUDDFFDDFF$,  we obtain the path $L'=U_3DUDDDFFFUDFFFFFFF$ (see Figure \ref{fig3bis} for an illustration of this example).

Now we will prove that any $U_kD$-equivalence class contains at most one element in $\mathcal{C}$.
For a contradiction, let $L$ and $L'$ be two different \L{}ukasiewicz paths in $\mathcal{C}$ belonging to the same class. We  write $L=QR$ and $L'=QS$ where $R$ and $S$ start with two different steps. Since $L$ and $L'$ lie in the same class, the two first steps of $R$ and $S$ cannot be $U_kD$ for $k\geq 1$. Moreover, since $L$ (resp. $L'$) lies into $\mathcal{C}$, the two first steps of $R$ (resp. $S$) cannot constituted a pattern $U_kF$ for $k\geq 1$. Then, $R$ and $S$ cannot start with any up step $U_k$, $k\geq 1$.

 Without loss of generality, let us assume that the first step of $R$ is a down step $D$ and then, the first step of $S$ is a flat step $F$.
 This means that the last point of $Q$ has its ordinate equal to zero (otherwise $L'$ could not belong to $\mathcal{C}$). As the first step of $R$ is $D$,
  the height of this step is $-1$  which gives a contradiction and completes the proof.\hfill $\Box$

\begin{thm}\label{thm7} The generating function for the set of  $U_kD$-equivalence classes of $\mathcal{L}$ with respect to the length is given by $$\frac{1-x +x^2+ \sqrt{1-2x-x^2-2x^3+x^4}}{1-2x-x^3+(1-x)\sqrt{1-2x-x^2-2x^3+x^4}},$$ which generates the generalized Catalan sequence defined by $g_0=1$, and   $g_{n+1} = g_n + \sum_{k=1}^{n-1} g_kg_{n-1-k}$ for $n\geq 0$
(see \href{https://oeis.org/A004148}{A004148} in \cite{Sloa}). \end{thm}
\noindent {\it Proof.} Using Lemma \ref{le3}, it suffices to obtain the generating function $C(x)$ for the set $\mathcal{C}$. A non-empty \L{}ukasiewicz path $L\in\mathcal{C}$
can be written either $L=FL'$ where $L'\in \mathcal{C}$, or $L=U_kDL_1 D L_2 D \ldots L_{k-1}DL'$ for $k\geq 1$ where
and $L_1, L_2, \ldots, L_{k-1}$ are some  \L{}ukasiewicz paths  without flats in $\bar{\mathcal{C}}$, and $L'\in \mathcal{C}$. So we obtain the functional equation
$C(x)=1+xC(x)+x^2C(x)\sum_{k\geq 0}x^k\bar{C}(x)^k$ where $\bar{C}(x)$ is the generating function for the set $\bar{\mathcal{C}}$ of \L{}ukasiewicz paths without flats in $\mathcal{C}$.
 Using the classical,
decomposition of a \L{}ukasiewicz path, we have $\bar{C}(x)=1+x\sum_{k\geq 1}x^k\bar{C}(x)^k$. A simple calculation provides the result.
\hfill $\Box$

\bigskip

\begin{thm}\label{thm7bis} For any $n\geq 0$, the map $\psi$ induces a bijection from $\mathcal{C}_n$ to the set of  Motzkin paths
 in $\mathcal{M}_n$  that avoid the pattern $UU$. \end{thm}
\noindent {\it Proof.} Theorem \ref{thm6bis}  ensures that $\psi$ is a bijection from $\mathcal{B}_n$ to $\mathcal{M}_n$ for $n\geq 0$. We have
$\mathcal{C}\subset\mathcal{B}$, and the paths in $\mathcal{C}$ are those in $\mathcal{B}$ that avoid the patterns $U_kU_\ell$, $k,\ell\geq 1$ and
$U_kF$, $k\geq 1$. Using Fact 2(b,d,e), the map $\psi$ transforms  occurrences of $U_kU_\ell$, $k,\ell\geq 1$, into occurrences of $UU$,
occurrences $U_kD$, $k\geq 2$, and  $U_kF$, $k\geq 1$ into occurrences of $UF$, and
occurrences of $UD$ into occurrences of $UD$. Then, the image by $\psi$ of $\mathcal{C}_n$ is
the subset of Motzkin paths in $\mathcal{M}_n$ that does not contain any pattern $UU$. \hfill $\Box$
\medskip

Let $\mathcal{E}\subset\mathcal{L}$ be the set of \L{}ukasiewicz paths such that any up step $U_k$, $k\geq 1$, is  immediately followed
 by a flat step $F$, and any flat step $F$ of positive height belongs to a pattern $U_kF$, $k\geq 1$.
For instance, we have $U_3DDDFUFD\notin\mathcal{E}$ and $U_3FDDFDUD\in\mathcal{E}$.
Let $\bar{\mathcal{E}}\subset\mathcal{E}$ be the set of \L{}ukasiewicz paths  in $\mathcal{E}$ without flat step on the $x$-axis.

\begin{lem}\label{le4}
There is a bijection between $\mathcal{E}$ and the set of $U_kF$-equivalence classes of $\mathcal{L}$.
\end{lem}

\noindent {\it Proof.} The proof is obtained {\it mutatis mutandis} as for Lemma \ref{le3} by replacing $U_kD$ with $U_kF$.
Let $L$ be a non-empty \L{}ukasiewicz path in $\mathcal{L}$. Let us prove that there exists a \L{}ukasiewicz path $L'\in\mathcal{E}$
(with the same length as $L$) such that $L$ and $L'$ belong to the same class.

We write
 $$L=K_0\prod_{i=1}^r\left(U_{k_i}FK_i\right),$$ where $k_i\geq 1$, $1\leq i\leq r$, and $K_0,K_1,K_2, \ldots, K_r$, $r\geq 0$,  are some parts (possibly empty)
 without pattern $U_kF$ for any $k\geq 1$.

We define the \L{}ukasiewicz path $$L'=F^{b_0}\prod_{i=1}^r\left(U_{k_1}FD^{a_i}F^{b_i-a_i}\right),$$ with $b_i=|K_i|$, $0\leq i\leq r$, and  for $1\leq i\leq r$,
 $a_i=\min\{b_i, k_i+\sum_{j=1}^{i-1} (k_j-a_j)\}$.
Less formally, $K_0$ is replaced with $F^{|K_0|}$, and $K_i$ is replaced with $D^{a_i}F^{b_i-a_i}$ where the value $a_i$, $1\leq i\leq r$, is the maximal number of down steps that can be placed between
 the two occurrences $U_{k_{i}}F$ and $U_{k_{i+1}}F$ so that $L'$ remains a lattice path.  Clearly, $L'\in\mathcal{E}$ and $L'$ belongs to the same class as $L$.

Now we will prove that any $U_kF$-equivalence class contains at most one element in $\mathcal{E}$.
For a contradiction, let $L$ and $L'$ be two different \L{}ukasiewicz paths in $\mathcal{E}$ belonging to the same class.
We  write $L=QR$ and $L'=QS$ where $R$ and $S$ start with two different steps.
Since $L$ and $L'$ lie in the same class, the two first steps of $R$ and $S$ cannot be $U_kF$ for $k\geq 1$.
So, since $L$ (resp. $L'$) lies into $\mathcal{E}$, any up step is followed by a flat step, which means that the
  first step of $R$ (resp. $S$) cannot be $U_k$, $k\geq 1$. Then, $R$ and $S$ cannot start with any up step $U_k$, $k\geq 1$.

 Without loss of generality, let us assume that the first step of $R$ is a down step $D$ and then, the first step of $S$ is a flat step $F$. Note that the first step
 of $S$ is necessarily on the $x$-axis.
 This means that the last point of $Q$ has its ordinate equal to zero (otherwise $L'$ could not belong to $\mathcal{E}$). As the first step
 of $R$ is $D$, the height of this step is $-1$  which gives a contradiction and completes the proof.
\hfill $\Box$

\begin{thm}\label{thm8} The generating function for the set of  $U_kF$-equivalence classes of $\mathcal{L}$ with respect to the length is given by
 $$\frac{1-x +2x^2+ \sqrt{1-2x+x^2-4x^3}}{1-2x-3x^3+(1-x+x^2)\sqrt{1-2x+x^2-4x^3}},$$ which generates the generalized Catalan sequence defined by $h_0=1$, and  for $n\geq 0$, $h_{n+1} = h_n + \sum_{k=0}^{n-2} h_kh_{n-2-k}$
  (see \href{https://oeis.org/A023431}{A023431} in \cite{Sloa}). \end{thm}

\noindent {\it Proof.} Using Lemma \ref{le4}, it suffices to obtain the generating function $E(x)$ for the set $\mathcal{E}$. A non-empty \L{}ukasiewicz path $L\in\mathcal{E}$ can be written
either $L=FL'$ where $L'\in \mathcal{E}$, or $L=U_kFL_1 D L_2 D \ldots L_{k}DL'$ for $k\geq 1$ where
and $L_1, L_2, \ldots, L_k$ are some  \L{}ukasiewicz paths in $\bar{\mathcal{E}}$, and $L'\in \mathcal{E}$. So we obtain the functional equation
$E(x)=1+xE(x)+x^3E(x)\sum_{k\geq 0}x^k\bar{E}(x)^{k+1}$ where $\bar{E}(x)$ is the generating function for the set $\bar{\mathcal{E}}$.
  Using the classical
decomposition of a \L{}ukasiewicz path, we have $\bar{E}(x)=1+\sum_{k\geq 3}x^k\bar{E}(x)^{k-1}$. A simple calculation provides the result.
\hfill $\Box$

\bigskip

\begin{thm}\label{thm8bis} For $n\geq 0$, the map $\psi$  is a bijection from $\mathcal{E}_n$ to the subset $\mathcal{M}'_n$ of
 Motzkin paths in $\mathcal{M}_n$  that avoid  $UU$ and $UD$. \end{thm}
\noindent {\it Proof.} We proceed by induction on $n$. Obviously, for  $n=0$, we have $\psi(\epsilon)=\epsilon$. For $0\leq k\leq n$, we assume that $\psi$ is a
bijection from $\mathcal{E}_k$ to the subset $\mathcal{M}'_k$  and we prove the result for $n+1$. Since the set of length $n$ Motzkin paths avoiding $UU$ and $UD$ is enumerated by the value $h_n$ defined in Theorem \ref{thm8}
(see \href{https://oeis.org/A023431}{A023431} in \cite{Sloa}), it suffices to prove that $\psi$ is surjective. So, let $M$ be a Motzkin path in  $\mathcal{M}'_{n+1}$. We distinguish two cases: ($i$) $M=FM'$ with $M'\in \mathcal{M}'_n$, and
 ($ii$) $M=UM'DM''$  where $M$ and $M'$ are two Motzkin paths in $\mathcal{M}'$.

 ($i$) Using the recurrence hypothesis, there is $L'\in \mathcal{E}_n$ such that $M'=\psi(L')$. So, the \L{}ukasiewicz path
 $L=FL'\in \mathcal{E}_{n+1}$ satisfies $\psi(L)=M$ which proves that $M$ belongs to the image by $\psi$ of $\mathcal{E}_{n+1}$.

 ($ii$) We suppose $M=UM'DM''$ with $M', M''\in\mathcal{M}'$. Since $M\in \mathcal{M}'_{n+1}$, we have $M'\neq \epsilon$ and $M'$ does not start with $U$, which implies that $M'$ starts with $F$. Using the recurrence
hypothesis, there are $L'\in \mathcal{E}$ and $L''\in \mathcal{E}$ such that $\psi(L')=M'$ and $\psi(L'')=M''$. Since $M'$ starts with a flat step, $L'$ also starts with a flat step.
 So,  $L=UL'DL''$ belongs to $\mathcal{E}_{n+1}$ and
satisfies $\psi(L)=M$ which proves that $\psi$ from $\mathcal{B}_n$ to $\mathcal{M}_n$ is bijective.
 \hfill $\Box$

\subsection{Modulo $FU_k$ for $k\geq 1$, and  $DU_k$ for $k\geq 1$  }
      For a given $\alpha\in \{D,F\}$, the $\alpha U_k$-equivalence in $\mathcal{L}$ is defined as follows: two \L{}ukasiewicz paths $L$ and $L'$
      of the same length are $\alpha U_k$-equivalent whenever for any $k\geq 1$, $L$ and $L'$  have the same  positions of $\alpha U_k$.

Let $\xi$ be the map from $\mathcal{L}$ to himself defined by $\xi(L)$ is obtained from $L$ by replacing any occurrence
 $U_kF$ by an occurrence $FU_k$ for $k\geq 1$. It is straightforward to verify that $\xi$ induces a bijection $\bar{\xi}$
  between the set of $U_kF$-equivalence classes and the set of $FU_k$-equivalence classes. Then, Theorem \ref{thm9} is directly deduced
  from Theorem \ref{thm8}.

\begin{thm}\label{thm9} The generating function for the set of  $FU_k$-equivalence classes of $\mathcal{L}$ with respect to the length also is the generating function given in Theorem \ref{thm8}. \end{thm}

Let $\mathcal{L}'_{n}$, $n\geq 2$, be the set of  \L{}ukasiewicz paths of length $n$ starting by $U$ and ending by $D$.
For $n\geq 0$, we define the bijection $\theta$ from $\mathcal{L}_n$ to $\mathcal{L}'_{n+2}$ as follows: $\theta(L)$ is obtained from $L$ by replacing
 any occurrence $U_kD$ by an occurrence $DU_k$ for $k\geq 1$, and by adding a step $U$ at the beginning and a step $D$ at the ending.
 It is straightforward to verify that $\theta$ induces a bijection $\bar{\theta}$ between the set of $U_kD$-equivalence classes of $\mathcal{L}_n$ and the set of $DU_k$-equivalence classes
 of $\mathcal{L}'_{n+2}$.
  Then, Theorem \ref{thm10} is deduced from Theorem \ref{thm7}.

\begin{thm}\label{thm10} The generating function for the set of  $DU_k$-equivalence classes of $\mathcal{L}$ with respect to the length
        is given by $$\frac{1-x -x^2-2x^3+ \sqrt{1-2x-x^2-2x^3+x^4}}{1-2x-x^3+(1-x)\sqrt{1-2x-x^2-2x^3+x^4}},$$ which generates the generalized Catalan sequence defined by
        $u_0=u_1=u_2=1$, and for $n\geq 3$  $u_{n} = g_{n-2}$ where $g_n$ is defined in Theorem \ref{thm7}.  (see \href{https://oeis.org/A004148}{A004148}  in \cite{Sloa}). \end{thm}

\subsection{Modulo $FF$}

Let $\mathcal{F}$ be the set constituted of  the union of $\{\epsilon, F\}$ with the set of \L{}ukasiewicz paths containing at most one up step $U_k$, $k\geq 1$ and such that any flat step $F$ is contained into a pattern $FF$ .
For instance,  $FFU_3DFFDDFFDFF\in\mathcal{F}$ and $FFU_3FFDDFDFFF\notin\mathcal{F}$.

\begin{lem}\label{le5}
There is a bijection between $\mathcal{F}$ and the set of $FF$-equivalence classes of $\mathcal{L}$.
\end{lem}

\noindent {\it Proof.} Let $L$ be a non-empty \L{}ukasiewicz path in $\mathcal{L}$. Let us prove that there exists a \L{}ukasiewicz path $L'\in\mathcal{F}$ (with the same length as $L$)
such that $L$ and $L'$ belong to the same class.
We write $$L=K_1F^{a_1}K_2F^{a_2}K_3\ldots K_rF^{a_r}K_{r+1},$$ with $r\geq 0$ and $a_i\geq 2$ for $1\leq i\leq r$, such that
 $K_1,K_2, \ldots, K_r,K_{r+1}$ are some parts  without pattern $FF$, $k\geq 1$, and $K_2,\ldots, K_r$ are not empty and do not have any $F$ in first and last position,
 and $K_1$ has no flat in last position, and $K_{r+1}$ has no flat in first position.

 If $L=F^n$ with $n\geq 0$, then its equivalence class is reduced to a singleton.
 Now let us assume that $L\neq F^n$. We distinguish two cases: (1) $K_1$ is not empty, and (2) $K_1$ is empty which means that $L=F^{a_1}K_2F^{a_2}K_3\ldots K_rF^{a_r}K_{r+1}.$

  In case (1), we define the \L{}ukasiewicz path $$L'=U_bD^{b_1-1}F^{a_1}D^{b_2}F^{a_2}D^{b_3}\ldots D^{b_r}F^{a_r}D^{b_{r+1}}$$  with $b_i=|K_i|$, $1\leq i\leq r+1$, and $b=b_1-1+\sum_{i=2}^{r+1}b_i$.

 In case (2), we define the \L{}ukasiewicz path $$L'=F^{a_1}U_bD^{b_2-1}F^{a_2}D^{b_3}\ldots D^{b_r}F^{a_r}D^{b_{r+1}}$$  with $b_i=|K_i|$, $2\leq i\leq r+1$, and $b=b_2-1+\sum_{i=3}^{r+1}b_i$.

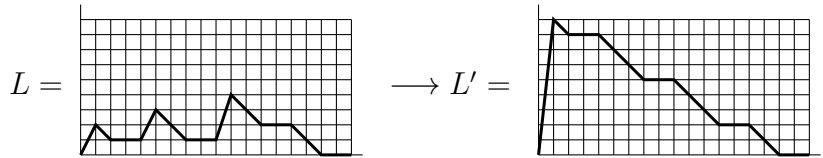
\begin{figure}[h]

 \begin{center} $L=$\begin{tabular}{c}\begin{tikzpicture}[scale=0.1]
            \draw (\A,\A)-- (38,\A);
             \draw (\A,\E)-- (\Ze,\E);
              \draw (\A,\C)-- (\Ze,\C);
               \draw (\A,\G)-- (\Ze,\G);
               \draw (\A,\I)-- (\Ze,\I);
              \draw (\A,\K)-- (\Ze,\K);
             \draw (\A,\M)-- (\Ze,\M);

             \draw (\A,\O)-- (\Ze,\O);\draw (\A,\Q)-- (\Ze,\Q);\draw (\A,\S)-- (\Ze,\S);

            \draw (\A,\A) -- (\A,\U);
             \draw (\C,\A) -- (\C,\S);\draw (\E,\A) -- (\E,\S);\draw (\G,\A) -- (\G,\S);
             \draw (\I,\A) -- (\I,\S);\draw (\K,\A) -- (\K,\S);\draw (\M,\A) -- (\M,\S);
             \draw (\O,\A) -- (\O,\S);\draw (\Q,\A) -- (\Q,\S);\draw (\S,\A) -- (\S,\S);
             \draw (\U,\A) -- (\U,\S);\draw (\W,\A) -- (\W,\S);\draw (\Y,\A) -- (\Y,\S);
             \draw (\ZZ,\A) -- (\ZZ,\S);
             \draw (\Za,\A) -- (\Za,\S);
             \draw (\Zb,\A) -- (\Zb,\S);
             \draw (\Zc,\A) -- (\Zc,\S);
             \draw (\Zd,\A) -- (\Zd,\S);
             \draw (\Ze,\A) -- (\Ze,\S);
            \draw[solid, line width=0.4mm] (\A,\A)--(\C,\E) -- (\E,\C)-- (\G,\C)  -- (\I,\C) -- (\K,\G)-- (\M,\E) -- (\O,\C) -- (\Q,\C)--(\S,\C) -- (\U,\I)-- (\W,\G) --(\Y,\E) -- (\ZZ,\E) --(\Za,\E) -- (\Zb,\C) -- (\Zc,\A) -- (\Zd,\A) -- (\Ze, \A);
         \end{tikzpicture}\\
         \end{tabular}
$\longrightarrow L'= $
          \begin{tabular}{c}\begin{tikzpicture}[scale=0.1]
            \draw (\A,\A)-- (38,\A);
             \draw (\A,\E)-- (\Ze,\E);
              \draw (\A,\C)-- (\Ze,\C);
               \draw (\A,\G)-- (\Ze,\G);
               \draw (\A,\I)-- (\Ze,\I);
              \draw (\A,\K)-- (\Ze,\K);
               \draw (\A,\M)-- (\Ze,\M);
            \draw (\A,\A) -- (\A,\U);
             \draw (\A,\O)-- (\Ze,\O);\draw (\A,\Q)-- (\Ze,\Q);\draw (\A,\S)-- (\Ze,\S);

             \draw (\C,\A) -- (\C,\S);\draw (\E,\A) -- (\E,\S);\draw (\G,\A) -- (\G,\S);
             \draw (\I,\A) -- (\I,\S);\draw (\K,\A) -- (\K,\S);\draw (\M,\A) -- (\M,\S);
             \draw (\O,\A) -- (\O,\S);\draw (\Q,\A) -- (\Q,\S);\draw (\S,\A) -- (\S,\S);
             \draw (\U,\A) -- (\U,\S);\draw (\W,\A) -- (\W,\S);\draw (\Y,\A) -- (\Y,\S);
             \draw (\ZZ,\A) -- (\ZZ,\S);
             \draw (\Za,\A) -- (\Za,\S);
             \draw (\Zb,\A) -- (\Zb,\S);
             \draw (\Zc,\A) -- (\Zc,\S);
             \draw (\Zd,\A) -- (\Zd,\S);
             \draw (\Ze,\A) -- (\Ze,\S);
 \draw[solid, line width=0.4mm] (\A,\A)--(\C,\S) -- (\E,\Q)-- (\G,\Q)  -- (\I,\Q) -- (\K,\O)-- (\M,\M) -- (\O,\K) -- (\Q,\K)--(\S,\K) -- (\U,\I)-- (\W,\G) --(\Y,\E) -- (\ZZ,\E) --(\Za,\E) -- (\Zb,\C) -- (\Zc,\A) -- (\Zd,\A) -- (\Ze, \A);
         \end{tikzpicture}
         \\
         \end{tabular}\\
\end{center}

         \caption{Illustration of the  example described in the proof of Lemma \ref{le5}.}
         \label{fig35}
\end{figure}

Less formally, we obtained $L'$ from $L$ by replacing any $K_i$ (excepted the first) with a run of down steps $D^{|K_i|}$, and by replacing the first $K_i$ ($K_1$ or $K_2$
according to the case (1) or (2)) with  $U_bD^{|K_1|-1}$ (or $U_bD^{|K_2|-1}$) where the up step $U_b$ balances all down steps in $L'$, {\it i.e.}, $b$ is the number of down steps in $L'$.
Clearly, $L'\in\mathcal{F}$ and $L'$ belongs to the same class as $L$. For instance, if  $L=U_2DFFU_2DDFFU_3DDFFDDFF$, then  $L'=U_9DFFDDDFFDDDFFDDFF$ (see Figure \ref{fig35}).

The definition of $\mathcal{L}$ implies that there is only one path of $\mathcal{F}$ in the same class as $L$, which completes the proof.\hfill $\Box$

\begin{thm}\label{thm11} The generating function for the set of  $FF$-equivalence classes of $\mathcal{L}$ with respect to the length is given by
 $${\frac {1-3x+4x^2-5x^3+7x^4-7x^5+6x^6-3x^7+{x}^{8}}{ \left(1-2x+x^2 -{x}^{3} \right)  \left( 1-
x \right) ^{2}}}.$$ (Note that the associated  sequence does not yet appear in \cite{Sloa}). \end{thm}

\noindent {\it Proof.} Using Lemma \ref{le5}, it suffices to obtain the generating function $F(x)$ for the set $\mathcal{F}$. A non-empty \L{}ukasiewicz path $L\in\mathcal{F}$ can be written either ($i$) $L=F^k$ for $k\geq 0$, or ($ii$) $L=F^{i_0} U_k F^{i_1}D^{j_1} F^{i_2}D^{j_2}\ldots F^{i_\ell}D^{j_\ell} F^{i_{\ell+1}}$  with $\ell\geq 1$, $i_0=0$ or $i_0\geq 2$, $i_1=0$ or $i_1\geq 2$, $i_{\ell+1}=0$ or $i_{\ell+1}\geq 2$, $i_m\geq 2$ for $2\leq m\leq \ell$, and $j_m\geq 1$ for $1\leq m\leq \ell$.

The generating function for the \L{}ukasiewicz paths satisfying ($i$) is given by $\frac{1}{1-x}$.

For \L{}ukasiewicz paths satisfying ($ii$), we give the generating function for each part of $L$, and we multiply them:

- For $F^{i_0}$, with $i_0=0$ or $i_0\geq 2$, the generating function  is $1+\frac{x^2}{1-x}$;

- For $F^{i_{\ell+1}}$, with $i_{\ell+1}=0$ or $i_{\ell+1}\geq 2$, the generating function is $1+\frac{x^2}{1-x}$;

- For $U_kF^{i_1}D^{j_1}$, with $i_1=0$ or $i_1\geq 2$ and $j_1\geq 1$, the generating function is $x(1+\frac{x^2}{1-x})\frac{x}{1-x}$;

- For  $F^{i_2}D^{j_2}\ldots F^{i_\ell}D^{j_\ell}$, with $i_m\geq 2$, and $j_m\geq 1$, the generating function is $\frac{1}{1-\frac{x^3}{(1-x)^2}}$.

Considering all these cases, we deduce:
 $$ F(x)= \left( 1+{\frac {{x}^{2}}{1-x}} \right) ^{3}{x}^{2} \left( 1-x
 \right) ^{-1} \left( 1-{\frac {{x}^{3}}{ \left( 1-x \right) ^{2}}}
 \right) ^{-1}+ \left( 1-x \right) ^{-1}$$ which completes the proof.\hfill $\Box$

\section{Concluding remarks}
Extending recent works on Dyck and Motzkin paths \cite{Bar1,Bar2},   the goal of this paper is to calculate the number of  \L{}ukasiewicz paths modulo the positions of a given pattern,
{\it i.e.} the number of possible sets $I=\{i_1,i_2, \ldots, i_k\}$
where $i_1,i_2,\ldots , i_k$  are the  occurrence positions of the pattern in \L{}ukasiewicz paths.
Can one do the same study for other lattice paths such as meanders, bridges and excursions, or Schroeder and Riordan paths?

For a pattern  $\alpha\in\{F,D,FD,DF,DD\}$, we have characterized the possible sets $I$ of positions of $\alpha$ in a \L{}ukasiewicz path.
More generally, is it possible  to  characterize these sets for  other patterns?   From our study, we can deduce a lower bound for the maximal cardinality of a class by calculating
  the average of cardinalities of the  classes, {\it i.e.},  the  total number of
 \L{}ukasiewicz paths divided by the number of classes.  Is it possible to calculate the exact value of  the maximal cardinality for a class, and for which set $I$ it is reached?
Also, it would be interesting to study some properties of the number of   \L{}ukasiewicz paths (of a
 given length) having $I$ as set of positions of the pattern. One can think  this number is a polynomial with respect to  the length $n$. If this is true, then we could give
  properties of these coefficients and roots, which would be a counterpart for lattice paths of the study of descent polynomial on the symmetric group $S_n$ (see MacMahon \cite{Mac}).


\end{document}